# On the Mannheim Surface Offsets

**Mehmet Önder[1], H. Hüseyin Uğurlu[2]**

[1] *Celal Bayar University, Faculty of Arts and Sciences, Department of Mathematics, 45047, Muradiye, Manisa, Turkey. E-mail:* mehmet.onder@cbu.edu.tr

[2] *Gazi University, Faculty of Education, Department of Secondary Education Science and Mathematics Teaching, Mathematics Teaching Program, Ankara, Turkey.*

## Abstract

In this paper, we study Mannheim surface offsets in dual space. By the aid of the E. Study Mapping, we consider ruled surfaces as dual unit spherical curves and define the Mannheim offsets of the ruled surfaces by means of dual geodesic trihedron (dual Darboux frame). We obtain the relationships between the invariants of Mannheim ruled surfaces. Furthermore, we give the conditions for these surface offset to be developable. Furthermore, we obtained that the dual spherical radius of curvature of offset surface is equal to dual offset angle.

## 1. Introduction

Generally, an offset surface is offset a specified distance from the original along the parent surface's normal. Offsetting of curves and surfaces is one of the most important geometric operations in CAD/CAM due to its immediate applications in geometric modeling, NC machining, and robot navigation [4]. Especially, the offsets of ruled surfaces, which are the surfaces generated by continuously moving of a straight line, have an important role in (CAGD) [11,12]. These surfaces are one of the most fascinating topics of surface theory and also used in many areas of science such as Computer Aided Geometric Design (CAGD), mathematical physics, moving geometry, kinematics for modeling the problems and model-based manufacturing of mechanical products. In [13], Ravani and Ku defined and given a generalization of the theory of Bertrand curve for Bertrand trajectory ruled surfaces on the line geometry. By considering the E. Study mapping, Küçük and Gürsoy have studied the integral invariants of closed Bertrand trajectory ruled surfaces in dual space [6]. They have given some characterizations of Bertrand offsets of trajectory ruled surfaces in terms of integral invariants (such as the angle of pitch and the pitch) of closed trajectory ruled surfaces and obtained the relationship between the area of projections of spherical images of Bertrand offsets of trajectory ruled surfaces and their integral invariants.

Recently, a new offset of ruled surfaces has been defined by Orbay, Kasap and Aydemir [7]. They have called this new offset as Mannheim offset and shown that every developable ruled surface has a Mannheim offset if and only if an equation should be satisfied between the geodesic curvature and the arc-length of spherical indicatrix of it. The corresponding characterizations of Mannheim offsets of ruled surfaces in Minkowski 3-space have been given in [8,9]. Furthermore, Mannheim offsets of ruled surfaces in dual space have been studied according to Blaschke frame in [10] and the characterizations of Mannheim offsets of ruled surfaces in terms of integral invariants of closed trajectory ruled surfaces have been given and the relationship between the area of projections of spherical images of Mannheim offsets of trajectory ruled surfaces and their integral invariants have been obtained.

In this paper, we examine the Mannheim offsets of the ruled surfaces in view of their dual geodesic trihedron (dual Darboux frame). Using the dual representations of the ruled surfaces, we give some theorems and new results characterizing Mannheim surface offsets.

## 3. Dual Representation of Ruled Surfaces

Dual numbers had been introduced by W.K. Clifford (1845-1879). A dual number has the form $\bar{a}=(a,a^*)=a+\varepsilon a^*$ where $a$ and $a^*$ are real numbers and $\varepsilon=(0,1)$ is dual unit with $\varepsilon^2=0$ [1]. The product of dual numbers $\bar{a}=(a,a^*)=a+\varepsilon a^*$ and $\bar{b}=(b,b^*)=b+\varepsilon b^*$ is given by

$$\bar{a}\bar{b}=(a,a^*)(b,b^*)=(ab,ab^*+a^*b)=ab+\varepsilon(ab^*+a^*b). \tag{1}$$

We denote the set of dual numbers by *D:*

$$D=\left\{\bar{a}=a+\varepsilon a^* : a,a^*\in\mathbb{R},\ \varepsilon^2=0\right\}. \tag{2}$$

Dual function of dual numbers presents a mapping of a dual numbers space on itself. Properties of dual functions were thoroughly investigated by Dimentberg [2]. He derived the general expression for dual analytic (differentiable) function as follows

$$f(\bar{x})=f(x+\varepsilon x^*)=f(x)+\varepsilon x^* f'(x), \tag{3}$$

where $f'(x)$ is derivative of $f(x)$. This definition allows us to write the dual forms of some well-known functions as follows

$$\begin{cases} \cos(\bar{x})=\cos(x+\varepsilon x^*)=\cos(x)-\varepsilon x^*\sin(x), \\ \sin(\bar{x})=\sin(x+\varepsilon x^*)=\sin(x)+\varepsilon x^*\cos(x), \\ \sqrt{\bar{x}}=\sqrt{x+\varepsilon x^*}=\sqrt{x}+\varepsilon\dfrac{x^*}{2\sqrt{x}},\ (x>0). \end{cases} \tag{4}$$

Let $D^3=D\times D\times D$ be the set of all triples of dual numbers, i.e.,

$$D^3=\left\{\tilde{a}=(\bar{a}_1,\bar{a}_2,\bar{a}_3):\quad \bar{a}_i\in D,\ i=1,2,3\right\}. \tag{5}$$

Then the set $D^3$ is called dual space. The elements of $D^3$ are called dual vectors. A dual vector $\tilde{a}$ may be expressed in the form $\tilde{a}=\vec{a}+\varepsilon\vec{a}^*=(\vec{a},\vec{a}^*)$, where $\vec{a}$ and $\vec{a}^*$ are the vectors of $\mathbb{R}^3$. Then for any vectors $\tilde{a}=\vec{a}+\varepsilon\vec{a}^*$ and $\tilde{b}=\vec{b}+\varepsilon\vec{b}^*$ in $D^3$, scalar product and cross product are defined by

$$\left\langle\tilde{a},\tilde{b}\right\rangle=\left\langle\vec{a},\vec{b}\right\rangle+\varepsilon\left(\left\langle\vec{a},\vec{b}^*\right\rangle+\left\langle\vec{a}^*,\vec{b}\right\rangle\right), \tag{6}$$

and

$$\tilde{a}\times\tilde{b}=\vec{a}\times\vec{b}+\varepsilon\left(\vec{a}\times\vec{b}^*+\vec{a}^*\times\vec{b}\right), \tag{7}$$

respectively, where $\left\langle\vec{a},\vec{b}\right\rangle$ and $\vec{a}\times\vec{b}$ are scalar product and cross product of the vectors $\vec{a}$ and $\vec{a}^*$ in $\mathbb{R}^3$, respectively.

The norm of a dual vector $\tilde{a}$ is given by

$$\|\tilde{a}\|=\sqrt{\langle\tilde{a},\tilde{a}\rangle}=\|\vec{a}\|+\varepsilon\frac{\left\langle\vec{a},\vec{a}^*\right\rangle}{\|\vec{a}\|}. \tag{8}$$

Then a dual vector $\tilde{a}$ with norm $1+\varepsilon 0$ is called dual unit vector and the set of dual unit vectors is

$$\tilde{S}^2=\left\{\tilde{a}=(\bar{a}_1,\bar{a}_2,\bar{a}_3)\in D^3 : \left\langle\tilde{a},\ \tilde{a}\right\rangle=1+\varepsilon 0\right\}, \tag{9}$$

which is called dual unit sphere [1,3].

In 3-dimensional space $\mathbb{R}^3$, an oriented line $L$ is determined by a point $p\in L$ and a unit vector $\vec{a}$. Then, one can define $\vec{a}^*=\vec{p}\times\vec{a}$ which is called moment vector. The value of $\vec{a}^*$ does not depend on the point $p$, because any other point $q$ in $L$ can be given by $\vec{q}=\vec{p}+\lambda\vec{a}$ and then $\vec{a}^* = \vec{p}\times\vec{a} = \vec{q}\times\vec{a}$. Reciprocally, given such a pair $(\vec{a},\vec{a}^*)$ one recovers the line $L$ as $L=\left\{(\vec{a}\times\vec{a}^*)+\lambda\vec{a} : \vec{a},\vec{a}^* \in \mathbb{R}^3,\ \lambda\in\mathbb{R}\right\}$, written in parametric equations. The vectors $\vec{a}$ and $\vec{a}^*$ are not independent of one another and they satisfy the following relationships

$$\langle\vec{a},\vec{a}\rangle = 1,\ \ \langle\vec{a},\vec{a}^*\rangle = 0. \tag{10}$$

The components $a_i,\ a_i^*\ (1\le i\le 3)$ of the vectors $\vec{a}$ and $\vec{a}^*$ are called the normalized Plucker coordinates of the line $L$. From (6), (9) and (10) we see that the dual unit vector $\tilde{a}=\vec{a}+\varepsilon\vec{a}^*$ corresponds to the line $L$. This correspondence is known as E. Study Mapping: There exists one-to-one correspondence between the vectors (or points) of dual unit sphere $\tilde{S}^2$ and the directed lines of the space $\mathbb{R}^3$ [1,3]. By the aid of this correspondence, the properties of the spatial motion of a line can be derived. Hence, the geometry of ruled surface is represented by the geometry of dual curves on the dual unit sphere in $D^3$.

The angle $\bar{\theta}=\theta+\varepsilon\theta^*$ between two dual unit vectors $\tilde{a},\tilde{b}$ is called dual angle and defined by

$$\langle\tilde{a},\tilde{b}\rangle = \cos\bar{\theta} = \cos\theta - \varepsilon\theta^*\sin\theta. \tag{11}$$

By considering the E. Study Mapping, the geometric interpretation of dual angle is that $\theta$ is the real angle between the lines $L_1, L_2$ corresponding to the dual unit vectors $\tilde{a},\ \tilde{b}$, respectively, and $\theta^*$ is the shortest distance between those lines [3].

In [14], Veldkamp introduced the dual geodesic trihedron of a ruled surface. Then, we speak to his procedure briefly as follows:

Let $(\tilde{k})$ be a dual curve represented by the dual vector $\tilde{e}(u)=\vec{e}(u)+\varepsilon\vec{e}^*(u)$. The unit vector $\vec{e}$ draws a curve on the real unit sphere $S^2$ and is called the (real) indicatrix of $(\tilde{k})$. We suppose throughout that it is not a single point. We take the parameter $u$ as the arc-length parameter $s$ of the real indicatrix and denote the differentiation with respect to $s$ by primes. Then we have $\langle\vec{e}',\vec{e}'\rangle=1$. The vector $\vec{e}'=\vec{t}$ is the unit vector parallel to the tangent of the indicatrix. The equation $\vec{e}^*(s)=\vec{p}(s)\times\vec{e}(s)$ has infinity of solutions for the function $\vec{p}(s)$. If we take $\vec{p}_o(s)$ as a solution the set of all solutions is given by $\vec{p}(s)=\vec{p}_o(s)+\lambda(s)\vec{e}(s)$, where $\lambda$ is a real scalar function of $s$. Therefore we have $\langle\vec{p}',\ \vec{e}'\rangle=\langle\vec{p}_o',\ \vec{e}'\rangle+\lambda$. By taking $\lambda=\lambda_o=-\langle\vec{p}_o',\ \vec{e}'\rangle$ we see that $\vec{p}_o(s)+\lambda_o(s)\vec{e}(s)=\vec{c}(s)$ is the unique solution for $\vec{p}(s)$ with $\langle\vec{c}',\vec{e}'\rangle=0$. Then, the given dual curve $(\tilde{k})$ corresponding to the ruled surface

$$\varphi_e=\vec{c}(s)+v\vec{e}(s), \tag{12}$$

may be represented by

$$\tilde{e}(s)=\vec{e}+\varepsilon\vec{c}\times\vec{e}, \tag{13}$$

where

$$\langle\vec{e},\vec{e}\rangle=1,\ \ \langle\vec{e}',\vec{e}'\rangle=1,\ \ \langle\vec{c}',\vec{e}'\rangle=0. \tag{14}$$

Then we have

$$\|\tilde{e}'\|=\vec{t}+\varepsilon\det(\vec{c}',\vec{e},\vec{t})=1+\varepsilon\Delta, \tag{15}$$

where $\Delta=\det(\vec{c}',\vec{e},\vec{t})$. The dual arc-length $\bar{s}$ of the dual curve $(\tilde{k})$ is given by

$$\bar{s}=\int_0^s\|\tilde{e}'(u)\|du=\int_0^s(1+\varepsilon\Delta)du=s+\varepsilon\int_0^s\Delta du. \tag{16}$$

Then $\bar{s}'=1+\varepsilon\Delta$. Therefore, the dual unit tangent to the curve $\tilde{e}(s)$ is given by

$$\frac{d\tilde{e}}{d\bar{s}}=\frac{\tilde{e}'}{\bar{s}'}=\frac{\tilde{e}'}{1+\varepsilon\Delta}=\tilde{t}=\vec{t}+\varepsilon(\vec{c}\times\vec{t})\,. \tag{17}$$

Introducing the dual unit vector $\tilde{g}=\vec{g}+\varepsilon\vec{c}\times\vec{g}$ we have the dual frame $\{\tilde{e},\tilde{t},\tilde{g}\}$ which is known as dual geodesic trihedron or dual Darboux frame of $\varphi_e$ (or $(\tilde{e})$). Also, it is well known that the real orthonormal frame $\{\vec{e},\vec{t},\vec{g}\}$ is called the geodesic trihedron of the indicatrix $\vec{e}(s)$ with the derivations

$$\vec{e}'=\vec{t},\;\; \vec{t}'=\gamma\vec{g}-\vec{e},\;\; \vec{g}'=-\gamma\vec{t}\,, \tag{18}$$

where $\gamma$ is called the conical curvature [5]. Similar to (18), the derivatives of the vectors of the dual frame $\{\tilde{e},\tilde{t},\tilde{g}\}$ are given by

$$\frac{d\tilde{e}}{d\bar{s}}=\tilde{t},\;\; \frac{d\tilde{t}}{d\bar{s}}=\bar{\gamma}\tilde{g}-\tilde{e},\;\; \frac{d\tilde{g}}{d\bar{s}}=-\bar{\gamma}\tilde{t}\,, \tag{19}$$

where

$$\bar{\gamma}=\gamma+\varepsilon(\delta-\gamma\Delta)\,,\quad \delta=\langle\vec{c}',\vec{e}\rangle, \tag{20}$$

and the dual darboux vector of the frame is $\tilde{d}=\bar{\gamma}\tilde{e}+\tilde{g}$. From the definition of $\Delta$ and (20), we also have

$$\vec{c}'=\delta\vec{e}+\Delta\vec{g}\,. \tag{21}$$

The dual curvature of dual curve(ruled surface) $\tilde{e}(s)$ is

$$\bar{R}=\frac{1}{\sqrt{(1+\bar{\gamma}^2)}}\,. \tag{22}$$

The unit vector $\tilde{d}_o$ with the same sense as the Darboux vector $\tilde{d}=\bar{\gamma}\tilde{e}+\tilde{g}$ is given by

$$\tilde{d}_o=\frac{\bar{\gamma}}{\sqrt{(1+\bar{\gamma}^2)}}\tilde{e}+\frac{1}{\sqrt{(1+\bar{\gamma}^2)}}\tilde{g}\,. \tag{23}$$

Then, the dual angle between $\tilde{d}_o$ and $\tilde{e}$ satisfies the followings

$$\cos\bar{\rho}=\frac{\bar{\gamma}}{\sqrt{(1+\bar{\gamma}^2)}},\;\; \sin\bar{\rho}=\frac{1}{\sqrt{(1+\bar{\gamma}^2)}}\,, \tag{24}$$

where $\bar{\rho}$ is the dual spherical radius of curvature. Hence $\bar{R}=\sin\bar{\rho},\; \bar{\gamma}=\cot\bar{\rho}$ (For details see [14]).

## 3. Characterizations of Mannheim Surface Offsets

Mannheim offsets of ruled surfaces have been defined by Orbay and et all as follows:

**Definition 3.1.** Assume that $\varphi$ and $\varphi^*$ be two ruled surfaces in $\mathbb{R}^3$ with the parametrizations

$$\varphi(s,v)=\vec{c}(s)+v\,\vec{q}(s),$$
$$\varphi^*(s,v)=\vec{c}^*(s)+v\,\vec{q}^*(s),$$

respectively, where $(\vec{c})$ (resp. $(\vec{c}^*)$) is the striction curve of the ruled surfaces $\varphi$ (resp. $\varphi^*$). Let Frenet frames of the ruled surfaces $\varphi$ and $\varphi^*$ be $\{\vec{q},\vec{h},\vec{a}\}$ and $\{\vec{q}^*,\vec{h}^*,\vec{a}^*\}$, respectively. The ruled surface $\varphi^*$ is said to be Mannheim offset of the ruled surface $\varphi$ if there exists a one to one correspondence between their rulings such that the asymptotic normal vector $\vec{a}$ of $\varphi$ is the central normal $\vec{h}^*$ of $\varphi^*$. In this case, $(\varphi,\varphi^*)$ is called a pair of Mannheim ruled surface [7].

Then, the dual version of Definition 3.1 can be given according to Darboux frame as follow:

**Definition 3.2.** Let consider the ruled surfaces $\varphi_e$ and $\varphi_{e_1}$ generated by dual unit vectors $\tilde{e}$ and $\tilde{e}_1$ and let $\{\tilde{e}(\bar{s}),\, \tilde{t}(\bar{s}),\, \tilde{g}(\bar{s})\}$ and $\{\tilde{e}_1(\bar{s}_1),\, \tilde{t}_1(\bar{s}_1),\, \tilde{g}_1(\bar{s}_1)\}$ be the dual Darboux frames of $\varphi_e$ and $\varphi_{e_1}$, respectively. Then $\varphi_e$ and $\varphi_{e_1}$ are called Mannheim surface offsets, if

$$\tilde{g}(\bar{s}) = \tilde{t}_1(\bar{s}_1)\,, \tag{25}$$

holds along the striction lines of the surfaces, where $\bar{s}$ and $\bar{s}_1$ are the dual arc-lengths of $\varphi_e$ and $\varphi_{e_1}$, respectively.

By this definition, the relationship between trihedrons of the ruled surfaces $\varphi_e$ and $\varphi_{e_1}$ can be given as follows

$$\begin{pmatrix} \tilde{e}_1 \\ \tilde{t}_1 \\ \tilde{g}_1 \end{pmatrix} = \begin{pmatrix} \cos\bar{\theta} & \sin\bar{\theta} & 0 \\ 0 & 0 & 1 \\ \sin\bar{\theta} & -\cos\bar{\theta} & 0 \end{pmatrix} \begin{pmatrix} \tilde{e} \\ \tilde{t} \\ \tilde{g} \end{pmatrix}, \tag{26}$$

where $\bar{\theta} = \theta + \varepsilon\theta^*$, $(0 \le \theta \le \pi,\ \theta^* \in \mathbb{R})$ is dual angle between the generators $\tilde{e}$ and $\tilde{e}_1$ of Mannheim ruled surface $\varphi_e$ and $\varphi_{e_1}$. The angle $\theta$ is called the offset angle and $\theta^*$ is called the offset distance. Then, $\bar{\theta} = \theta + \varepsilon\theta^*$ is called dual offset angle of the Mannheim ruled surface $\varphi_e$ and $\varphi_{e_1}$. If $\theta = 0$ and $\theta = \pi/2$ then the Mannheim offsets are said to be oriented offsets and right offsets, respectively [10].

Now, we give some theorems and results characterizing Mannheim surface offsets.

**Theorem 3.1.** *Let $\varphi_e$ and $\varphi_{e_1}$ form a Mannheim surface offset. The offset angle $\theta$ and the offset distance $\theta^*$ are given by*

$$\theta = -s + c,\quad \theta^* = -\int_0^s \Delta du + c^*, \tag{27}$$

respectively, where $c$ and $c^*$ are real constants.

**Proof.** Suppose that $\varphi_e$ and $\varphi_{e_1}$ form a Mannheim offset. Then from (26) we have

$$\tilde{e}_1 = \cos\bar{\theta}\,\tilde{e} + \sin\bar{\theta}\,\tilde{t}\ . \tag{28}$$

Differentiating (28) with respect to $\bar{s}$ gives

$$\frac{d\tilde{e}_1}{d\bar{s}} = -\sin\bar{\theta}\left(1 + \frac{d\bar{\theta}}{d\bar{s}}\right)\tilde{e} + \cos\bar{\theta}\left(1 + \frac{d\bar{\theta}}{d\bar{s}}\right)\tilde{t} + \bar{\gamma}\sin\bar{\theta}\,\tilde{g}\ . \tag{29}$$

Since $\dfrac{d\tilde{e}_1}{d\bar{s}}$ and $\tilde{g}$ are linearly dependent, from (29) we get $\dfrac{d\bar{\theta}}{d\bar{s}} = -1$. Then for the dual constant $\bar{c} = c + \varepsilon c^*$ we write

$$d\bar{\theta} = -d\bar{s},$$

$$\bar{\theta} = -\bar{s} + \bar{c},$$

$$\theta + \varepsilon\theta^* = -s - \varepsilon s^* + c + \varepsilon c^*,$$

and from (16) we have

$$\theta = -s + c,\ \theta^* = -\int_0^s \Delta du + c^*,$$

where $c$ and $c^*$ are real constants.

**Corollary 3.1.** *Let $\varphi_e$ and $\varphi_{e_1}$ form a Mannheim surface offset. Then $\varphi_e$ is developable if and only if $\theta^* = c^* = \text{constant}$.*

**Proof.** Since $\varphi_e$ and $\varphi_{e_1}$ form a Mannheim surface offset we have Theorem 4.1. Thus from (27) we see that $\varphi_e$ is developable i.e. $\Delta = 0$ if and only if $\theta^* = c^* = \text{constant}$.

**Theorem 3.2.** *Let $\varphi_e$ and $\varphi_{e_1}$ form a Mannheim surface offset. Then there is the following differential relationship between the dual arc-length parameters of $\varphi_e$ and $\varphi_{e_1}$*

$$\frac{d\bar{s}_1}{d\bar{s}} = \bar{\gamma}\sin\bar{\theta}. \tag{30}$$

**Proof.** Suppose that $\varphi_e$ and $\varphi_{e_1}$ form a Mannheim offset. Then, from Theorem 3.1 we have

$$\frac{d\tilde{e}_1}{d\bar{s}_1} = \tilde{t}_1 = \bar{\gamma}\sin\bar{\theta}\frac{d\bar{s}}{d\bar{s}_1}\tilde{g}. \tag{31}$$

From (26) we have $\tilde{t}_1 = \tilde{g}$. Then (31) gives us

$$\bar{\gamma}\sin\bar{\theta}\frac{d\bar{s}}{d\bar{s}_1} = 1, \tag{32}$$

and from (32) we get (30).

**Theorem 3.3.** *Let $\varphi_e$ and $\varphi_{e_1}$ form a Mannheim surface offset. Then there are the following relationships between the real arc-length parameters and invariants of $\varphi_e$ and $\varphi_{e_1}$*

$$\begin{cases} \dfrac{ds_1}{ds} = \gamma\sin\theta, \\ \Delta_1 = \theta^*\cot\theta + \dfrac{\delta}{\gamma}. \end{cases} \tag{33}$$

**Proof.** Let $\varphi_e$ and $\varphi_{e_1}$ form a Mannheim surface offset. Then from Theorem 3.2, (30) holds. By considering (20) the real and dual parts of (30) are

$$\frac{ds_1}{ds} = \gamma\sin\theta,\ \frac{ds\,ds_1^* - ds^*ds_1}{ds^2} = \theta^*\gamma\cos\theta + (\delta - \gamma\Delta)\sin\theta, \tag{34}$$

respectively. Furthermore from (16) we have

$$ds^* = \Delta ds,\ ds_1^* = \Delta_1 ds_1. \tag{35}$$

Writing the equalities (35) in (34) we have (33).

**Corollary 3.2.** *Let $\varphi_e$ and $\varphi_{e_1}$ form a Mannheim surface offset. Then, the Mannheim offset $\varphi_{e_1}$ is developable if and only if*

$$\theta^* = -\frac{\delta}{\gamma}\tan\theta, \tag{36}$$

*holds.*

**Proof.** Assume that $\varphi_e$ and $\varphi_{e_1}$ form a Mannheim surface offset. From Theorem 3.3, the condition to satisfy $\Delta_1 = 0$ is that

$$\theta^* = -\frac{\delta}{\gamma}\tan\theta. \tag{37}$$

**Theorem 3.4.** *Let $\varphi_e$ and $\varphi_{e_1}$ form a Mannheim surface offset. Then*

$$\delta_1 = \frac{\delta}{\gamma}\cot\theta - \theta^* \tag{38}$$

*holds.*

**Proof.** Let the striction lines of $\varphi_e$ and $\varphi_{e_1}$ be $c(s)$ and $c_1(s_1)$, respectively and let $\varphi_e$ and $\varphi_{e_1}$ form a Mannheim surface offset. Then we can write

$$\vec{c}_1 = \vec{c} + \theta^* \vec{g}\,. \tag{39}$$

Differentiating (39) with respect to $s_1$ we have

$$\frac{d\vec{c}_1}{ds_1} = \left(\frac{d\vec{c}}{ds} - \theta^*\gamma\vec{t} + \frac{d\theta^*}{ds}\vec{g}\right)\frac{ds}{ds_1}. \tag{40}$$

From (20) we know that $\delta_1 = \langle d\vec{c}_1 / ds_1, e_1\rangle$. Then from (26) and (40) we obtain

$$\delta_1 = \left(\cos\theta\langle d\vec{c}/ds,\, \vec{e}\rangle + \sin\theta\langle d\vec{c}/ds,\, \vec{t}\rangle - \theta^*\gamma\sin\theta\langle\, \vec{t},\, \vec{t}\,\rangle\right)\frac{ds}{ds_1}. \tag{41}$$

Since $\delta = \langle d\vec{c}/ds, e\rangle$ and $\langle d\vec{c}/ds,\, \vec{t}\,\rangle = 0$, from (41) we write

$$\delta_1 = \left(\delta\cos\theta - \theta^*\gamma\sin\theta\right)\frac{ds}{ds_1}. \tag{42}$$

Furthermore, from (33) we have

$$\frac{ds}{ds_1} = \frac{1}{\gamma\sin\theta}, \tag{43}$$

and substituting (43) in (42) we obtain

$$\delta_1 = \frac{\delta}{\gamma}\cot\theta - \theta^*.$$

**Theorem 3.5.** *Let $\varphi_e$ and $\varphi_{e_1}$ form a Mannheim surface offset. The conical curvature of $\varphi_{e_1}$ is obtained as follows,*

$$\gamma_1 = \cot\theta. \tag{44}$$

**Proof.** From (18) and (26) we have

$$\begin{aligned}\gamma_1 &= -\langle \vec{g}_1', \vec{t}_1\rangle \\ &= -\left\langle \frac{d}{ds_1}(\sin\theta\vec{e} - \cos\theta\vec{t}),\, \vec{g}\right\rangle\end{aligned} \tag{45}$$

which gives

$$\gamma_1 = \gamma\cos\theta\frac{ds}{ds_1}. \tag{46}$$

From the first equality of (33) and (46) we have (44).

**Theorem 3.6.** *Let* $\varphi_e$ *and* $\varphi_{e_1}$ *form a Mannheim surface offset. Then, the dual curvature of* $\varphi_{e_1}$ *is given by*

$$\overline{R}_1 = \sin\overline{\theta}\,. \tag{47}$$

**Proof.** From (22) and (44), Eq. (47) is obtained immediately.

**Theorem 3.7.** *Let* $\varphi_e$ *and* $\varphi_{e_1}$ *form a Mannheim surface offset. Then, the dual spherical radius of curvature of* $\varphi_{e_1}$ *is equal to dual offset angle, i.e.,*

$$\overline{\rho}_1 = \overline{\theta}\,. \tag{48}$$

**Proof.** From (24) we have $\sin\overline{\rho}_1 = \overline{R}_1$. Then by (4) and (47) we get

$$\sin\rho_1 = \sin\theta, \quad \rho_1^* \cos\rho_1 = \theta^* \cos\theta\,. \tag{49}$$

Then we have

$$\rho_1 = \theta, \;\; \rho_1^* = \theta^*\,. \tag{50}$$

which gives $\overline{\rho}_1 = \overline{\theta}$.

From Eq. (23) and Theorem 3.6, we have the following corollary.

***Corollary 3.3.*** *Let* $\varphi_e$ *and* $\varphi_{e_1}$ *form a Mannheim surface offset. Then the relationship between dual offset angle* $\overline{\theta}$ *and dual unit darboux vector* $\tilde{d}_{0_1}$ *of* $\varphi_{e_1}$ *is given as follows*

$$\tilde{d}_{0_1} = \cos\overline{\theta}\,\tilde{e}_1 + \sin\overline{\theta}\,\tilde{g}_1\,.$$

**Example 3.1.** Let consider the hyperbolic paraboloid surface $\varphi_e$ given by the parametrization

$$\varphi_e(s,v) = \left(\frac{1}{2}s, \frac{1}{2}s, 0\right) + v\left(\frac{1}{2}, -\frac{1}{2}, s\right), \tag{51}$$

and rendered in Fig. 1.

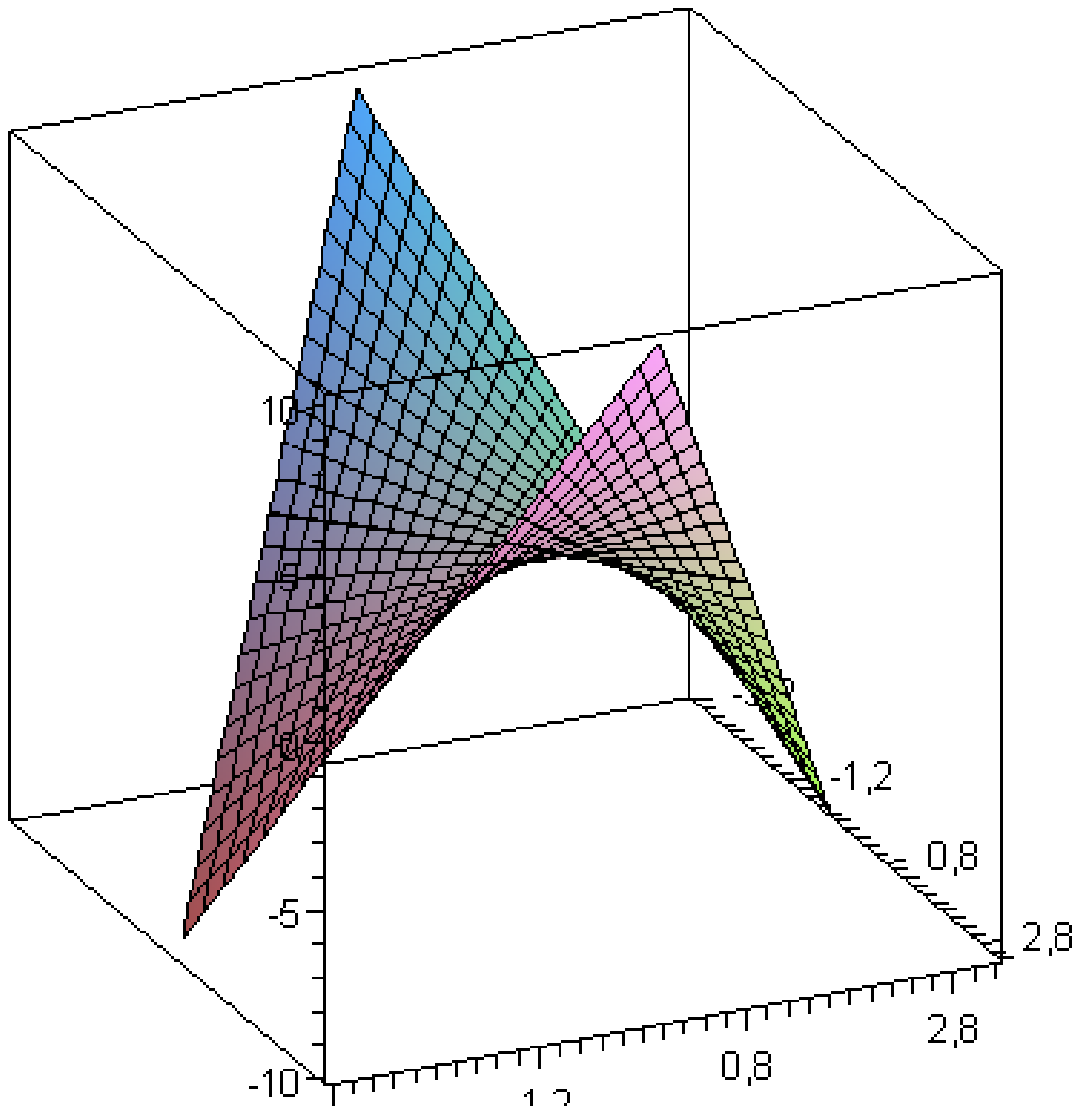

**Fig. 1** Hyperbolic paraboloid surface $\varphi_e$.

From the E. Study Mapping, the dual spherical curve representing (51) is

$$\tilde{e}(s)=\frac{\sqrt{2}}{\sqrt{1+2s^2}}\left[\left(\frac{1}{2},-\frac{1}{2},s\right)+\varepsilon\left(\frac{1}{2}s^2,-\frac{1}{2}s^2,-\frac{1}{2}s\right)\right]. \tag{52}$$

Then, the dual Darboux frame of $\varphi_e$ is obtain as follows

$$\begin{cases}\tilde{e}(s)=\dfrac{\sqrt{2}}{\sqrt{1+2s^2}}\left[\left(\dfrac{1}{2},-\dfrac{1}{2},s\right)+\varepsilon\left(\dfrac{1}{2}s^2,-\dfrac{1}{2}s^2,-\dfrac{1}{2}s\right)\right]\\ \tilde{t}(s)=\dfrac{1}{\sqrt{1+8\tan^2(\sqrt{2}s)}}\left[\left(-2\tan(\sqrt{2}s),2\tan(\sqrt{2}s),1\right)+\varepsilon\left(\tan(\sqrt{2}s),-\tan(\sqrt{2}s),4\tan^2(\sqrt{2}s)\right)\right]\\ \tilde{g}(s)=\left(-\dfrac{\sqrt{2}}{2},-\dfrac{\sqrt{2}}{2},0\right)\end{cases}$$

The general equation of the Mannheim offset surface of $\varphi_e$ is

$$\begin{aligned}\varphi_{e_1}(s,v)=&\left(\frac{1}{2}s-\theta^*\frac{\sqrt{2}}{2},\frac{1}{2}s-\theta^*\frac{\sqrt{2}}{2},0\right)\\ &+v\left(\frac{\sqrt{2}}{\sqrt{1+2s^2}}\cos\theta\left(\frac{1}{2},-\frac{1}{2},s\right)+\sin\theta\left(-\frac{\sqrt{2}}{2},-\frac{\sqrt{2}}{2},0\right)\right)\end{aligned} \tag{53}$$

From (53) we can give the following special cases:

**i)** The Mannheim offset $\varphi_{e_1}$ with dual offset angle $\bar{\theta}=0+\varepsilon 4\sqrt{2}$ is

$$\varphi_{e_1}(s,v)=\left(\frac{1}{2}s-4,\frac{1}{2}s-4,0\right)+v\left(\frac{\sqrt{2}}{\sqrt{1+2s^2}},-\frac{\sqrt{2}}{\sqrt{1+2s^2}},\frac{\sqrt{2}s}{\sqrt{1+2s^2}}\right)$$

which is an oriented offset of $\varphi_e$ (Fig. 2).

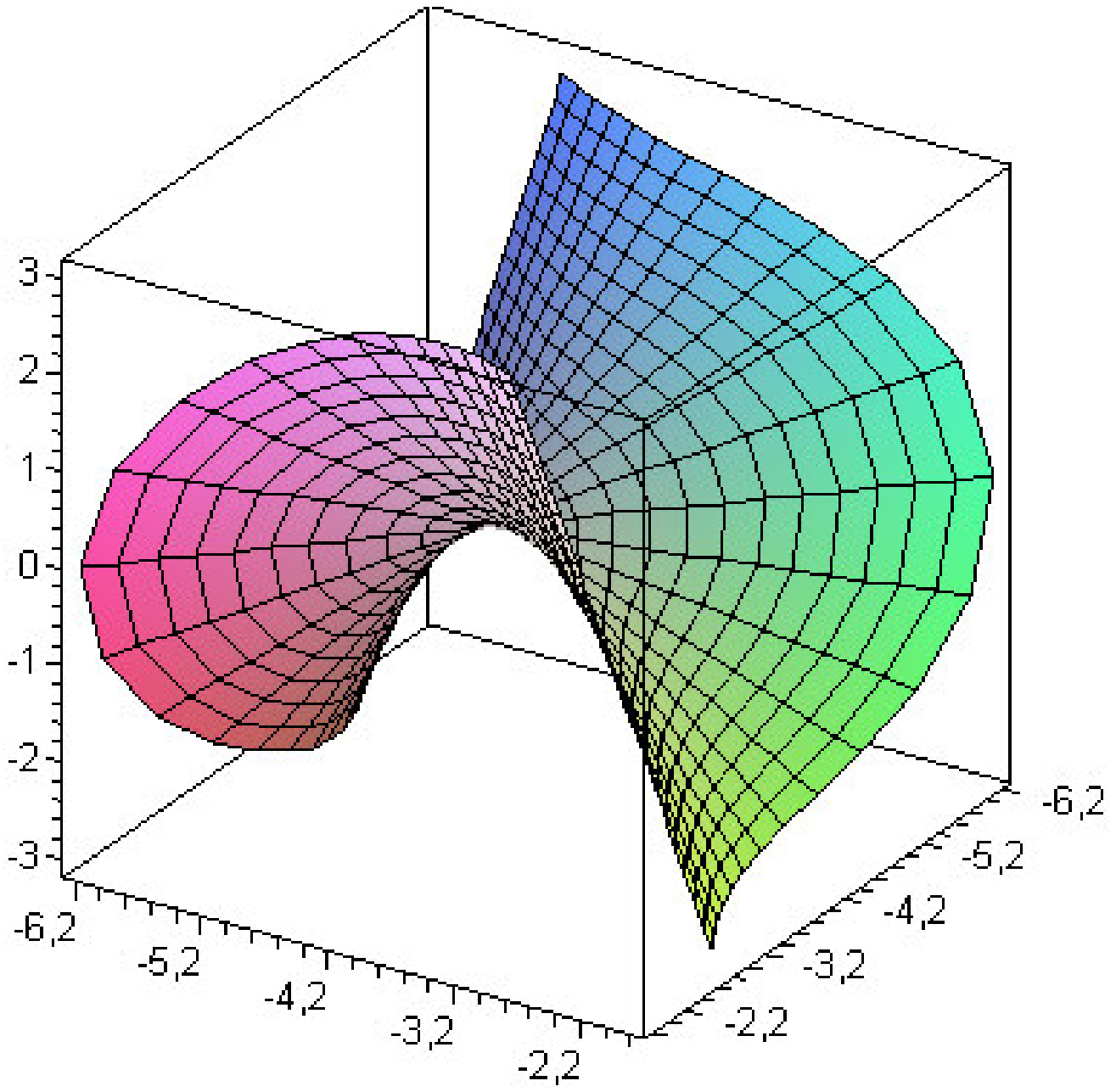

**Fig. 2.** Mannheim offset $\varphi_{e_1}$ with dual offset angle $\bar{\theta}=0+\varepsilon 4\sqrt{2}$

**ii)** The Mannheim offset $\varphi_{e_1}$ with dual offset angle $\bar{\theta} = \pi/4 + \varepsilon 2\sqrt{2}$ is

$$\varphi_{e_1}(s,v) = \left(\frac{1}{2}s - 2, \frac{1}{2}s - 2, 0\right) + v\left(\frac{1}{2\sqrt{1+2s^2}} - \frac{1}{2}, \frac{-1}{2\sqrt{1+2s^2}} - \frac{1}{2}, \frac{s}{\sqrt{1+2s^2}}\right)$$

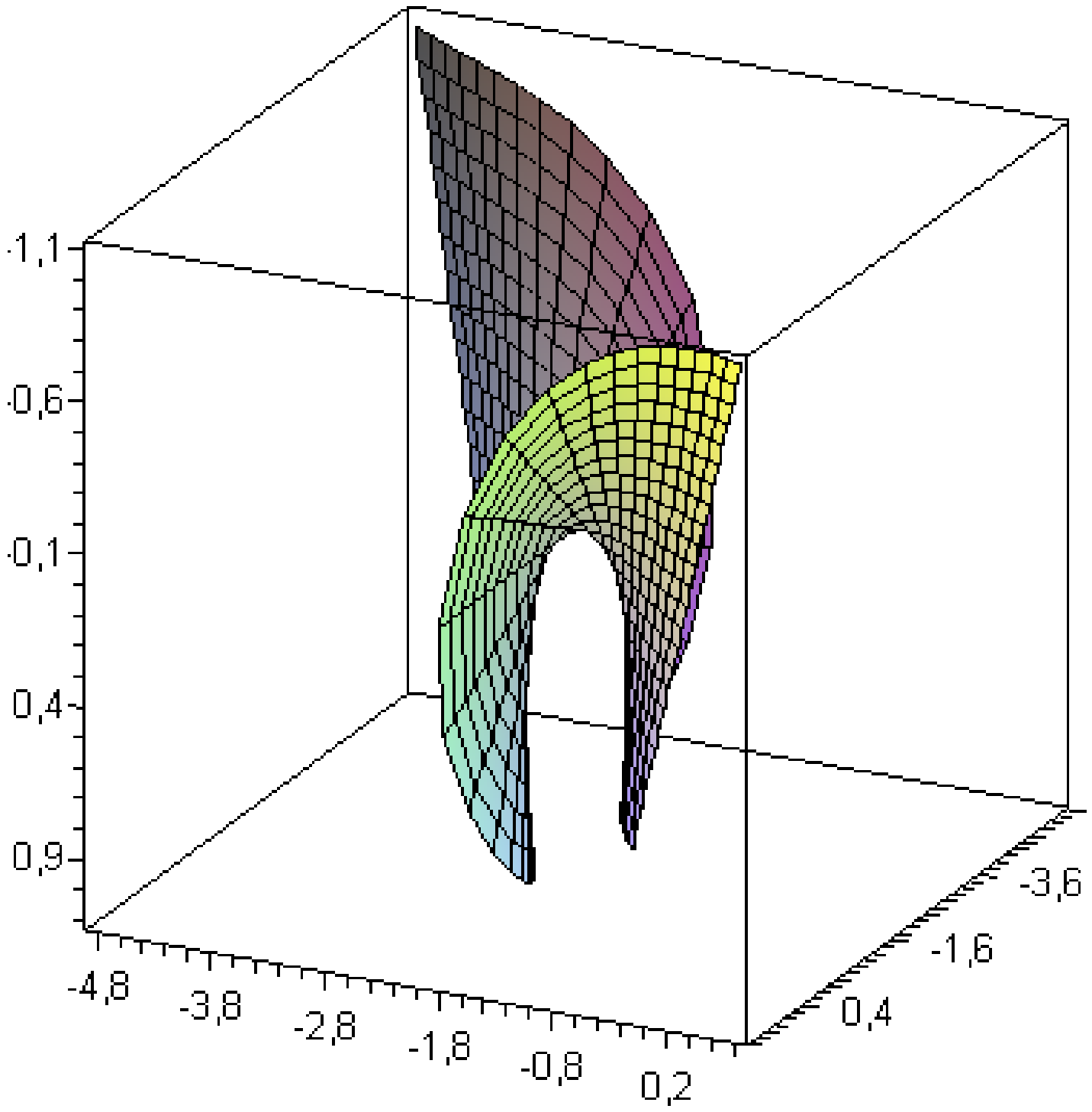

**Fig. 3.** Mannheim offset $\varphi_{e_1}$ with dual offset angle $\bar{\theta} = \pi/4 + \varepsilon 2\sqrt{2}$

## 4. Conclusions

In this paper, we give the characterizations of Mannheim offsets of ruled surfaces in view of dual geodesic trihedron (dual Darboux frame). We find new relations between the invariants of Mannheim surface offsets. Furthermore, we give the relationships for Mannheim surface offsets to be developable.